\documentclass[11pt]{amsart}
\usepackage{amssymb,amsmath,amsthm}
\usepackage{epsfig}
\newtheorem{thm}{Theorem}
\newtheorem{lemma}[thm]{Lemma}

\newtheorem{defn}[thm]{Definition}
\newtheorem{conj}[thm]{Conjecture}
\newtheorem{prop}[thm]{Proposition}
\newtheorem{cor}[thm]{Corollary}

\newtheorem{ex}[thm]{Example}
\newtheorem{formula}[thm]{}

\newcommand\Sk{{\mathcal S}}

\title{Skein Modules of 3-Manifolds}
\author{J\'ozef H. Przytycki}
\begin{document}
\maketitle
\centerline{Presented by A. Bia{\l}ynicki-Birula on January 25, 1989}
\ \\
\ \\
\begin{quotation}
Summary. \baselineskip=10pt

It is natural to try to place the new polynomial invariants of links in
algebraic topology (e.g. to try to interpret them using homology or
homotopy groups). However, one can think that these new polynomial
invariants are byproducts of a new more delicate algebraic invariant of
3-manifolds which measures the obstruction to isotopy of links (which
are homotopic). We propose such an algebraic invariant based on skein
theory introduced by Conway (1969) and developed by Giller (1982) as
well as Lickorish and Millett (1987).
\\ \
\end{quotation}
\ \\
Let $M$ be an oriented 3-manifold and $R$ a commutative ring with $1$.
For $r_0,\ldots,r_{k-1} \in R$ we define the $k$th skein module
${\mathcal S}_k(M;R)(r_0,\ldots,r_{k-1})$ as follows:\\
Let ${\mathcal L}(M)$ be the set of all ambient isotopy classes of
oriented links in $M$. Let ${\mathcal M}({\mathcal L},R) $ be a
free $R$-module generated by ${\mathcal L}(M)$ and
${\mathcal S}_{{\mathcal L}(M)}(r_0,\ldots,r_{k-1})$
the submodule generated by linear skein expressions
$r_0 L_O+r_1 L_1+ \ldots + r_{k-1}L_{k-1}$,
 where $ L_0,L_1, \ldots ,L_{k-1}$ are classes of links identical except the
 parts shown in Fig. 1.\\
\ \\
\centerline{\psfig{figure=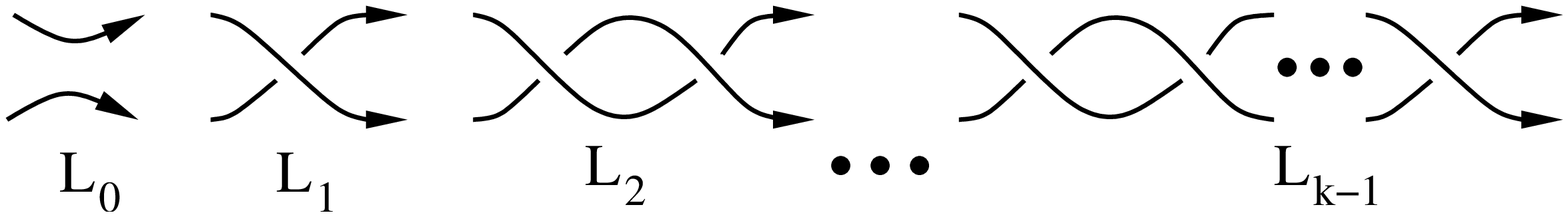,height=1.7cm}}
\begin{center}
Fig. 1
\end{center}

\begin{defn}\label{Definition 1}
 The $R$-module $${\mathcal S}_k (M;R)(r_0,\ldots,r_{k-1})
= {\mathcal M}({\mathcal L},R)/S_{{\mathcal
L}(M)}(r_0,\ldots,r_{k-1})$$ is called the kth skein module of $M$.

\end{defn}
\begin{ex}\label{Example 2}\
\begin{itemize}
  \item [(a)] ${\mathcal S}_k (M;R)(0,\ldots,0)
= {\mathcal M}({\mathcal L},R)$,
  \item [(b)] ${\mathcal S}_k (M;R)(1,0,-1)$ is a free $R$ module with
  basis consisting of links up to homotopy (i.e. two links
  $\alpha_1,\alpha_2: S^1\cup S^1\cup \ldots \cup S^1
  \rightarrow M$ are equivalent if they are homotopic).
\end{itemize}
\end{ex}
\begin{thm}\label{Theorem 3}
$\Sk_3(S^3;R)(r_0,r_1,r_2)=R$ provided that $r_0$, $r_1$ and $r_2$ are
invertible in R.
\end{thm}
Theorem 3 is the reformulation of the existence of the skein
polynomial [see 2,10].\\
\begin{prop}\label{Proposition 4}
We have the following elementary properties of skein modules:
\begin{itemize}
  \item [(a)] If $f:M\rightarrow N$ is an embedding of 3-manifolds, then $f$
  induces the homomorphism of skein modules:
  $$f_*: \Sk_k(M;R)(r_0,\ldots,r_{k-1})\rightarrow S_k(N;R)
(r_0,\ldots,r_{k-1}).$$
  \item [(b)] If $f, g:M\rightarrow N$ are isotopic embeddings, then $f_*=g_*$.
  \item [(c)] If $N$ is obtained from a 3-manifold M by adding a 2-handle
  along a simple closed curve in $\partial M$ and $f:M\rightarrow N$
  is the natural embedding then
   $f_*: \Sk_k(M; R)(r_0,\ldots,r_{k-1})\rightarrow
  \Sk_k(N; R)(r_0,\ldots,r_{k-1})$ is an epimorphism.
  \item [(d)] If $\widehat{M}$ is obtained from $M$ by cupping
  off 2-spheres in $
  \partial M$ by 3-cells and $f:M\rightarrow \widehat{M}$ is the natural
  embedding then
  $$f_*: \Sk_k(M; R)(r_0,\ldots,r_{k-1})\rightarrow
  \Sk_k(\widehat{M};R)(r_0,\ldots,r_{k-1})$$
  is an isomorphism.
  \item [(e)] For a disjoint sum of 3-manifolds $M \sqcup N$ one has:
  $$ \Sk_k(M \sqcup N; R)(r_0,\ldots,r_{k-1})=
  \Sk_k(M; R)(r_0,\ldots,r_{k-1})\otimes \Sk_k(N; R)(r_0,\ldots,r_{k-1}).$$
\end{itemize}
\end{prop}
One also has the following universal coefficients lemma:
\begin{lemma}\label{Lemma 5}
\begin{itemize}
  \item [(a)] Let $\varphi: R \rightarrow R^{'}$ be a homomorphism of rings
 and $r_0,\ldots,r_{k-1}$ elements of $R$. Then
  $$ {\mathcal S}_k(M; R^{'})(\varphi(r_0),\ldots,\varphi(r_{k-1}))=
 {\mathcal S}_k(M; R)(r_0,\ldots,r_{k-1})\otimes_R R^{'}$$
 where R isomorphism (as well as $R^{'}$ isomorphism) between modules
 is defined by identity on classes of links in M.\\
In particular:\ \
  \item [(b)] $\Sk_k(M; R)(r_0,\ldots,r_{k-1})=$\\
  $ {\mathcal S}_k(M; Z[x_0,\ldots,x_{k-1}])(x_0,\ldots,x_{k-1})
   \otimes_{Z[x_0,\ldots,x_{k-1}]}R.$
  \item [(c)] If $r_0,\ldots,r_{k-1}$ are invertible in $R$, then:
${\mathcal S}_k(M;
  R)(r_0,\ldots,r_{k-1})=
  {\mathcal S}_k(M;Z[x_0^{\mp 1},\ldots,x_{k-1}^{\mp 1}])(x_0,\ldots,x_{k-1})
  \otimes_{Z[x_0^{\mp 1},\ldots,x_{k-1}^{\mp 1}]}R$.\\
  In (b), (respectively (c)), ${\mathcal S}_k(M; R)(r_0,\ldots,r_{k-1})$ and
  $R$ have the structure of
  $Z[x_0,\ldots,x_{k-1}]$ (resp. $Z[x_0^{\mp 1},\ldots,x_{k-1}^{\mp
  1}]$) module due to the homomorphism given by $\varphi (x_i)=r_i$.
\end{itemize}
\end{lemma}
It follows from Lemma 5(b) and Example 2(a) that\\
 $ {\mathcal S}_k(M;Z[x_0,\ldots,x_{k-1}])(x_0,\ldots,x_{k-1})$
fully classifies
links in M. That is, non-isotopic links represent different
elements in the skein module.

The most convenient choice of $R$,
which allows some computations but is still quite general, is
$Z[x_0^{\mp 1},\ldots,x_{k-1}^{\mp
  1}]$. We will denote the skein module\\
  ${\mathcal S}_k(M;Z[x_0^{\mp 1},\ldots,x_{k-1}^{\mp 1}])(x_0,\ldots,x_{k-1})$
  shortly by ${\mathcal S}_k(M)$.

If we assume that $M$ is compact-oriented 3-manifold with no
2-spheres in the boundary and $k>2$, then ${\mathcal S}_k(M)$  is 
known only for $k=3, M={\mathcal S}^3$ (Theorem 3)
and $k=3,\ M={\mathcal S}^1\times D^2$.
This latter
can be derived from the result obtained in \cite{Ho-K}.
\begin{thm}(see \cite{Ho-K})\label{Theorem 6}
${\mathcal S}_3 (S^1 \times D^2)$ is a free $Z[x_0^{\mp
1},\ldots,x_{k-1}^{\mp 1}]$ module with basis consisting of
a trivial circle and families of layered torus links of type $(k,1)$
$k\neq 0$. The layered family of (k,1)-torus links
$\{H_{i_1},H_{i_2},\ldots, H_{i_n}\}$ is constructed as follows:

Let $S^1 \times D^2=S^1 \times  \langle  0,1\rangle \times \langle
0,2n \rangle  $ and $H_{i_j}$
 be a torus knot of type $(i_j,1)$ on the torus
  $\partial (S^1 \times  \langle  0,1 \rangle \times  
\langle 2j-1,2j \rangle )$ then
  $\{H_{i_1},\ldots, H_{i_n}\}$ is equal to
$H_{i_1} \cup \ldots \cup H_{i_n}$ and we disregard the order of
$H_{i_j}$'s (compare Fig. 2).

\end{thm}
\ \\
\centerline{\psfig{figure=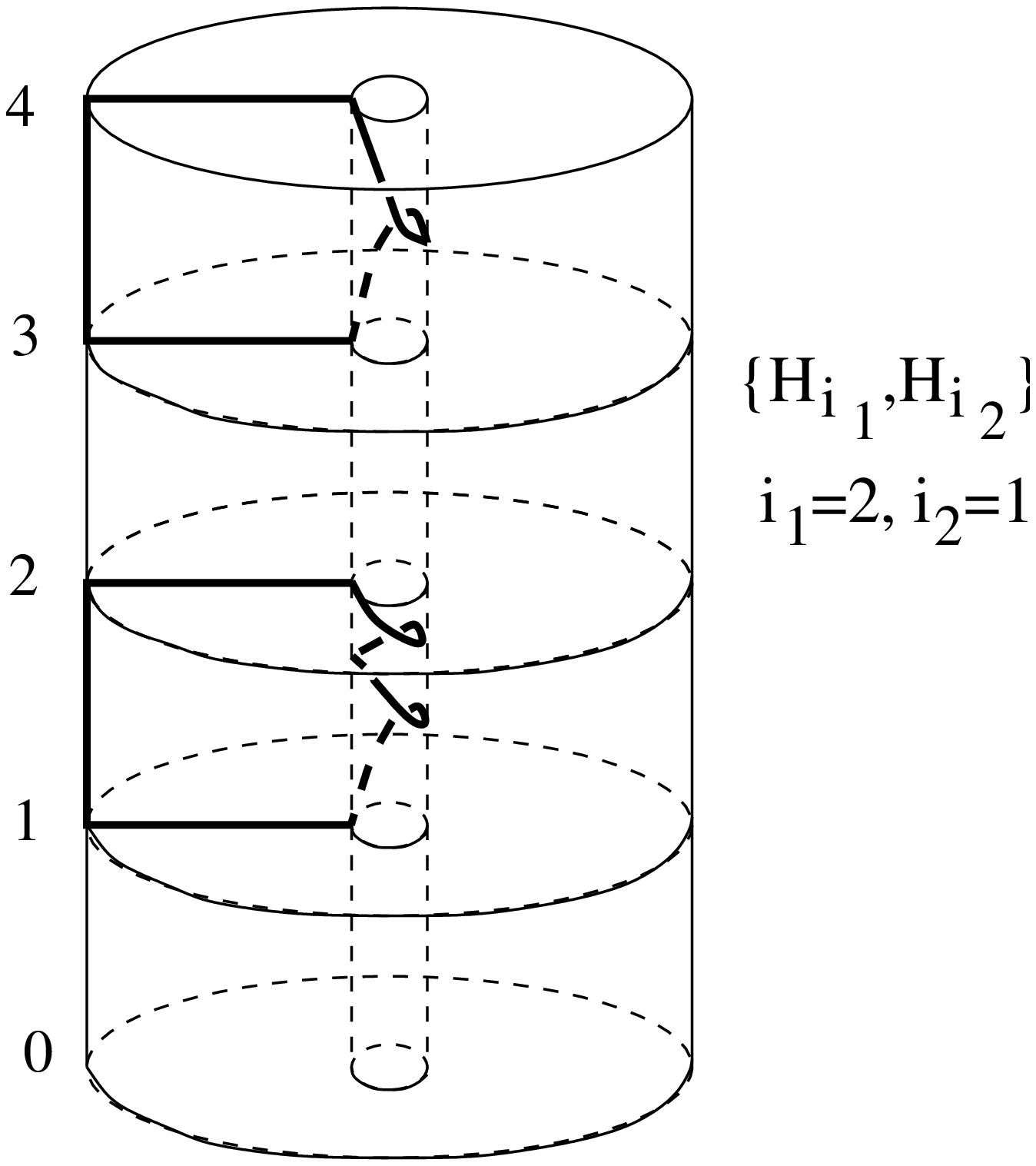,height=8.6cm}}
\centerline{$S^1 \times D^2 = S^1 \times  \langle 0,1 \rangle
\times  \langle 0,4 \rangle $}
\begin{center}
Fig. 2
\end{center}

\begin{conj}\label{Conjecture 7}
Let $ \langle  G\rangle$ denote the set of conjugacy classes of
elements of $G$ different than 1. Then for a handlebody $H_n$,
 ${\mathcal S}_3(H_n)$ is the free
$Z[x_0^{\mp 1},x_1^{\mp 1},x_2^{\mp 1}]$-module with basis
$(\cup_{i=1}^{\infty} \mathrm{X}_i  \langle \pi_1 (H_n) \rangle )/\sim $
where $\sim $ identifies sequences which differ only by
permutation. It can also be described as a module of polynomials
$Z[x_0^{\mp 1},x_1^{\mp 1},x_2^{\mp 1}][ \langle \pi_1(H_n) \rangle ]$,
with variables from $ \langle \pi_1(H_n) \rangle $, as suggested
by A. Joyal. The empty sequence (or polynomial 1) corresponds to a
trivial knot\footnote{Added for e-print:\
 I didn't utilize fully the suggestion of A.~Joyal,
otherwise I would realize that the empty sequence should correspond
to the empty link and that the skein module for the product of
a surface and the interval is not only the ring (as discussed
before Figure 5) but also a unital algebra with the empty link as
the identity element. The lack of the identity element made my description
of the structure of the Kauffman bracket skein module (and ring) of
a solid torus unnecessary awkward.}. A sequence
$\gamma_1, \ldots, \gamma_s$ (or a monomial $\gamma_1 \ldots
\gamma_s$) corresponds to a layered family of links representing
$\gamma_1 \ldots \gamma_s$ in $H_n$ decomposed into product of a
planar surface and an interval. It should be stressed that (unlike
in $H_1$) we do not expect $\gamma_i, \gamma_j$ to represent
always the same element of the skein module as $\gamma_j,
\gamma_i$. Therefore, different choice of sequence representatives
may lead to different basis of the skein module.
\end{conj}

This conjecture is a generalization of Theorems 3 and 6. 

It is not true, 
however, that if $\pi_1 (M)$ is free then ${\mathcal S}_3(M)$ is free.
 $M=S^1 \times S^2$ is the example. The knot $K=S^1 \times *$ 
represents a nontrivial
 element in ${\mathcal S}_3(S^1 \times S^2)$
but $(x_0+x_1+x_2)(x_0-x_1+x_2)K$ is equal to $0$
 in ${\mathcal S}_3(S^1 \times S^2)$.
We can further conjecture that ${\mathcal S}_3(M)$ is always a
 quotient module of $Z[x_0,x_1,x_2][ \langle \pi_1(M) \rangle ]$ where as in Conjecture 7 the
 monomials should correspond to layered families of links in a handlebody
being a part of a Heegaard decomposition of $M$.

The next conjecture is a weak version of the Van Kampen-Seifert
Theorem. Namely, if $M^\Delta N$ is a disk sum of two 3-manifolds
with connected boundary then the natural embedding $i:M
\rightarrow M^\Delta N$ induces
a monomorphism on skein modules.\\

{\bf Relative skein modules.} \ \
Let $F$ be a 2-submanifold of a  3-manifold $M$ where $F$ is either in
the boundary of $M$ or $F$ is properly embedded in $M$ (that is
$\partial M \cap F=\partial F$). Similarly as before, $R$ is
 a commutative ring with $1$ and
 $r_0, \ldots, r_{k-1} \in R$. We define the kth relative skein module
${\mathcal S}_k(M,F;R)(r_0,\ldots,r_{k-1})$ as follows. Let $\mathcal
{L}(M,F)$ denote the set of all ambient isotopy classes $\mod F$ of
oriented compact 1-manifolds $L$ in $M$ such that $L \cap F= \partial
F$. Elements of $\mathcal L(M,F)$ are called relative links in
$(M,F)$. Now $\mathcal M (\mathcal L(M,F); R)$ is a free $R$-module
generated by ${\mathcal L}(M,F)$ and
${\mathcal S}_{{\mathcal L}(M,F)}(r_0,\ldots,r_{k-1})$ is its submodule
generated by linear skein expression
\begin{formula}\label{Formula 8}
\ \ \ \ \ \ \ \ \ \ \ \ \ \ \ \ \ \ \ 
$ r_0L_0+ r_1L_1+ \ldots + r_{k-1}L_{k-1}$
\end{formula}
where diagrams of relative links $L_0, \ldots, L_{k-1}$ are
considered to be different in the exterior\footnote{Added for
e-print: they can be different only in the exterior of $F$,
that is, they are identical on $F$.} of $F$.
\begin{defn}\label{Definition 9}
The module $${\mathcal S}_k(M, F; R)(r_0,\ldots, r_{k-1})= \mathcal M
(\mathcal L(M,F);R) / S_{\mathcal L(M,F)}(r_0,r_1,\ldots,
r_{k-1})$$ is called the kth relative skein module of $(M,F)$.
\end{defn}
The points $L \cap F$ are called inputs and outputs of L. The
expression (8) uses relative links with the same set of inputs,
say $a_1, \ldots, a_n$, and outputs, say $b_1, \ldots, b_n$.
Therefore, for the given set of inputs and outputs one has the
well-defined submodule ${\mathcal S}_k(M,F;a_1, \ldots, a_n,b_1,\ldots,
b_n;R)(r_0,\ldots,r_{k-1})$ of the module
${\mathcal S}_k(M,F;R)(x_0,\ldots,x_{k-1})$ generated by relative links with
$a_1,\ldots, a_n$ as inputs and $b_1,\ldots, b_n$ as outputs. This
concept was first considered in \cite{L-M} (preskein, linearization of a
room, etc.)
\begin{ex}\label{Example 10}
Let $k=3$, $F=Z[a^{\mp 1},z^{\mp 1}]$, $r_0=a^{-1},r_1=-z,r_2=a$ (i.e.
the relation is a $L_+ +a^{-1}L_-=zL_0$), $M=D^3$, $F=\partial D^3$,
and $F$ has $n$ inputs $(a_1,\ldots, a_n)$ and $n$ outputs
$(b_1, \ldots, b_n)$ then
$${\mathcal S}_3(D^3, \partial D^3, a_1,\ldots, a_n, b_1,\ldots,
b_n; Z[a^{\mp 1},z^{\mp 1}])(a^{-1},-z, a) $$
is a free $Z[a^{\mp 1}, z^{\mp 1}]$ module of rank $n!$. One can construct
a base of the
module by choosing representatives of $n!$ pairing inputs and outputs.\\
For a proof see \cite{M-T}. 
In the notation of \cite{L-M} we have computed the
linearization of a room with n inputs and n outputs (Fig. 3).
\end{ex}
\ \\
\centerline{\psfig{figure=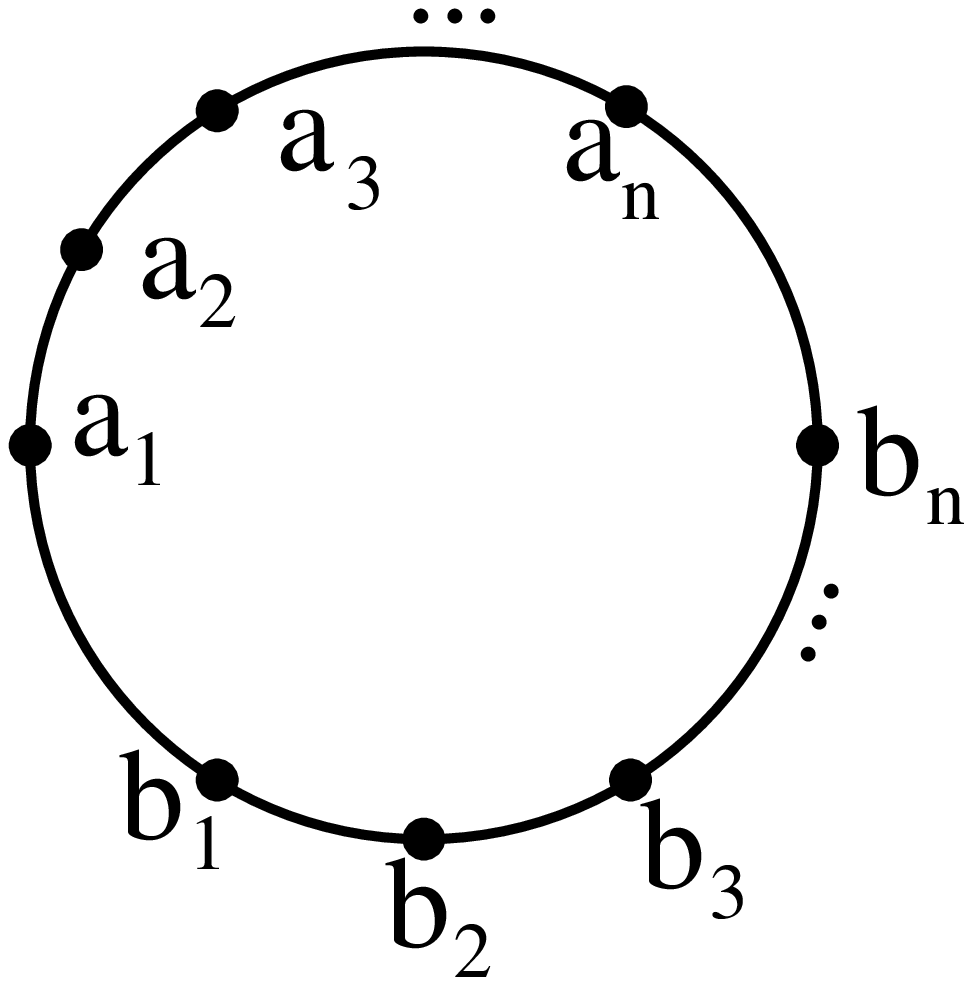,height=5.8cm}}
\begin{center}
Fig. 3
\end{center}

One can modify the notion of skein modules by considering diagrams:\\

\ \\
\centerline{\psfig{figure=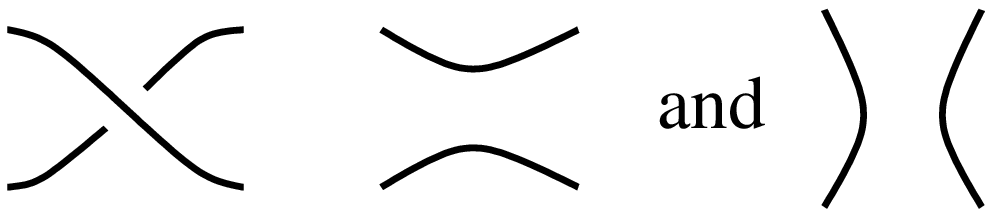,height=1.6cm}}
\vspace*{0.1cm}
\ \\
and the Kauffman approach. We will define skein module
${\mathcal S}_{2,\infty}(M,V_1;R)(A)$ of the Heegaard decomposition $(M,
V_1,\overline{M-V_1})$ of $M$ where $V_1$ is a handlebody and $M$
is obtained
from $V_1$ by adding 2 (and 3) handles to $\partial V_1$.\\

We start the definition from the handlebody $H_n$.

\begin{defn}\label{Definition 11}
Let $H_n=F_{0,n+1}\times [-1,1]$ where $F_{0,n+1}$ is a disk with
n holes and let $p: H_n \rightarrow F_{0,n+1}$ be the natural
projection. Let R be a ring as in Definition 1 and $\mathcal
M(\mathcal L_N,R)$ a free R module generated by $\mathcal L_N(H_n)$-
the set of all ambient isotopy classes of unoriented links in
$H_n$. Define ${\mathcal S}_{2,\infty}(A)$ for $A$ invertible in $R$,
to be the submodule of $\mathcal M(\mathcal L_n,R)$ generated by the
expression:
 $(-A^3)^{swL_+}L_+ - (-A^3)^{swL_0} A L_0 -
(-A^3)^{swL_{\infty}} A^{-1}L_{\infty}$
 where $L_+,L_0,L_{\infty}$ are classes of unoriented links represented by
 diagrams on $F_{0,n+1}$ identical except for parts shown in Fig.4;
$sw(L)$ is a self-writhe of L i.e.
\ $sw(L)=\sum_q sgn\ q$\  where the sum is taken over all
 selfcrosssings of $L$ ($sw(L)=w(L)-2lk(L)$ for any orientation of $L$). 
Now define\footnote{Added for e-print:\ My definition is rather cumbersome 
as I work with unframed diagrams. Notice that $(-A^3)^{swL}L$ 
can be viewed as framed link with ``blackboard" framing.}
 ${\mathcal S}_{2,\infty}(H_n,R)(A)=
\mathcal M(\mathcal L_N,R)/S_{2,\infty}(A)$.
\end{defn}

For any Heegaard decomposition $(M,V_1,\overline {M-V_1})$ of $M$, one
can define ${\mathcal S}_{2,\infty}(M,V_1;R)(A)$ by dividing
${\mathcal S}_{2,\infty}(V_1,R)(A)$ by relations on
 $\mathcal M(\mathcal L,R)$ which identify two elements
of $\mathcal L_N(V_1)$
which are isotopic in $M$.\\

\ \\
\centerline{\psfig{figure=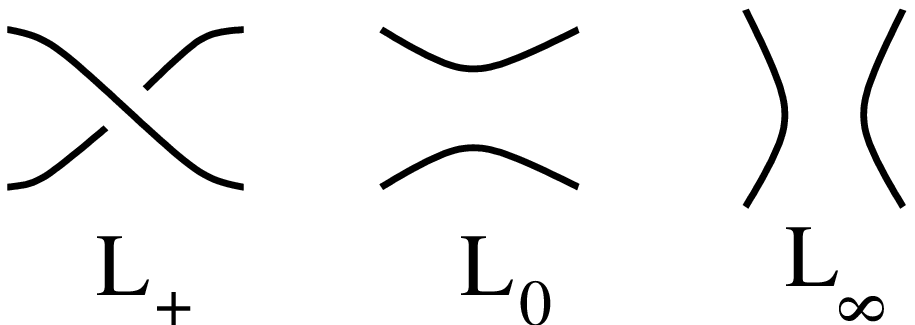,height=2.2cm}}
\begin{center}
Fig. 4
\end{center}

It is difficult to compute the skein module in general, but for
handlebodies it is an easier task.

\begin{thm}\label{Theorem 11}
\begin{itemize}
\item [(a)] ${\mathcal S}_{2,\infty}(S^3,R)(A)=R$
\item [(b)] ${\mathcal S}_{2,\infty} (H_n,R)(A)$
is a free R module with free basis $\mathcal L_N(F_{0,n+1})$
consisting of a trivial knot in $F_{0,n+1}$ and links in
$F_{0,n+1}$ with no trivial components.
\end{itemize}
\end{thm}
Proof:
\begin{itemize}
\item[(a)] can be thought as a special case of (b). In fact, the
definition of
  ${\mathcal S}_{2,\infty}(M,R)$ was chosen so that the Kauffman proof
of the existence of
  the bracket polynomial \cite{Kauf} gives us Theorem 12(a).
More precisely, if
  $R=Z[A^ {\mp 1}]$ and the trivial circle is the generator of
${\mathcal S}_{2,\infty}(S^3,R)$ then the coefficient in $R$ of
any link represented by a diagram $L$ is equal to $(-A^3)^{-swL}
 \langle L \rangle $ where $ \langle  L \rangle $ is the Kauffman
bracket of $L$.
\end{itemize}
The proof of (b) proceeds in two steps:\\

\begin{itemize}
    \item [Step 1.] Two diagrams of links in $F_{0,n+1}$ represent
(ambient) isotopic links if they are equivalent by a sequence of
Reidemeister moves of type $\Omega_1^{\mp 1}, \Omega_2^{\mp 1}$ or
$\Omega_3^{\mp 1}$. It is essentially the Reidemeister Theorem \cite{Rei}.
    \item [Step 2.] Now we follow Kauffman approach. $L_+$ is a
    linear combination of $L_0$ and $L_{\infty}$, therefore
    ${\mathcal L}_N(F_{0,n+1})$
     generates ${\mathcal S}_{2,\infty}(H_n,R)(A)$. On the other hand one can
     easily check (see [6]) that each link diagram on $F_{0,n+1}$
     up to Reidemeister moves can be uniquely expressed as an
     $R$-linear combination of elements of $\mathcal
     L_N(F_{0,n+1})$, therefore $\mathcal L(F_{0,n+1})$ forms a
     basis for ${\mathcal S}_{2,\infty}(H_n,R)(A)$.\\

\end{itemize}

Now we will analyze more precisely the situation of a  solid torus
$(H_1=S^1 \times D^2= F_{0,2} \times [-1,1])$ (cf.\cite{H-P}). The
basis ${\mathcal L}_N(F_{0,2})$ of
${\mathcal S}_{2,\infty}(S^1 \times D^2,R)(A)$ consists
of a trivial knot in $F_{0,2}$ and copies of circles parallel to
$\partial F_{0,2}$ (Fig. 5).

The choice of basis depends upon the product structure
$$ S^1 \times D^2= F_{0,2} \times [-1,1]$$

If one considers a link in $S^3$ one of whose components, say $K$,
is a trivial circle, then $K$ determines the product structure of
the solid torus $S^1 \times D^2= S^3- intV_k$, where $V_k$ is a
tubular neighborhood of $K$ in $S^3$ (the meridian of $V_k$ is a
longitude of $S^1 \times D^2$). Let $L$ be the rest of the link.
Thus, $L \subset S^1 \times D^2$.
\ \\
\centerline{\psfig{figure=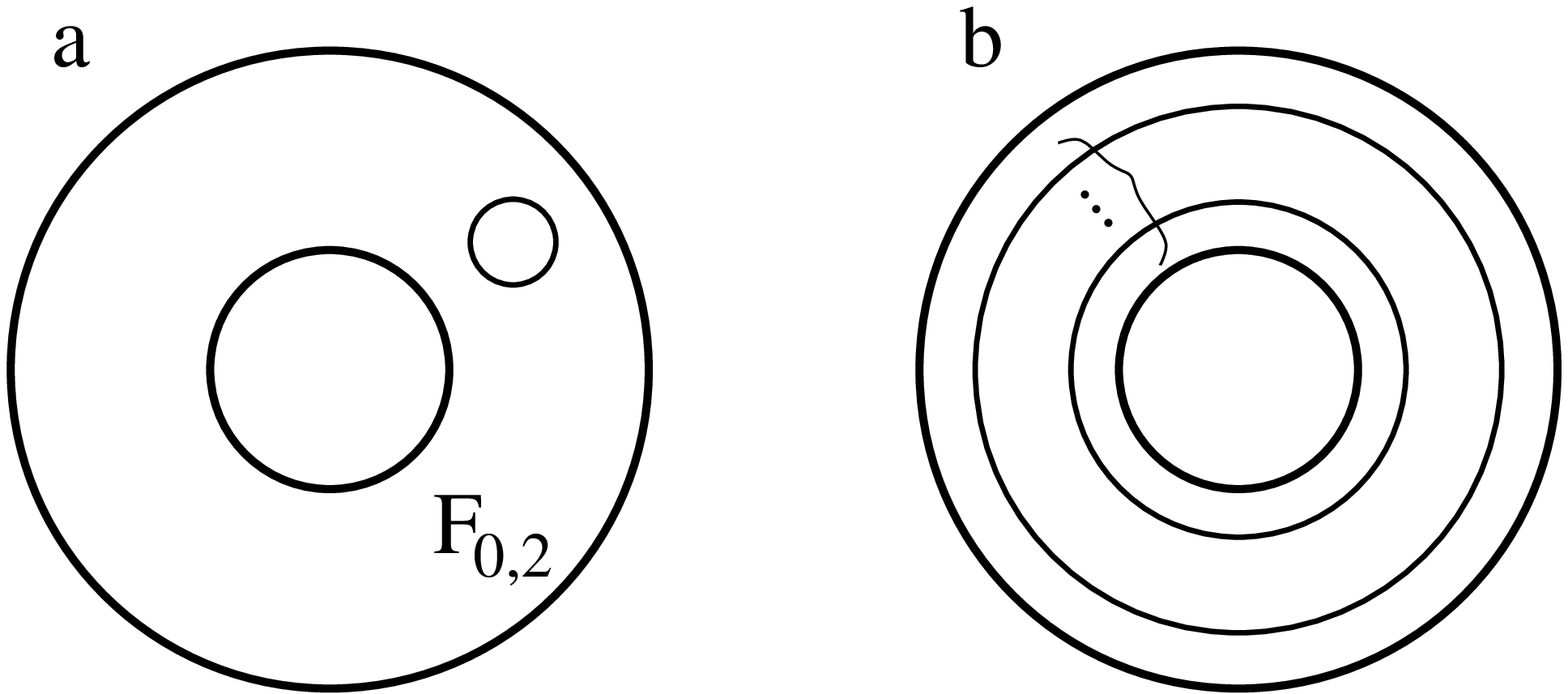,height=3.7cm}}
\begin{center}
Fig. 5
\end{center}

The elements of basis $\mathcal L(F_{0,2})$ can be denoted by
$h^0=1, h, -h^2(A^2+A^{-2}), \ldots,h^i(-A^2-A^{-2})^{i-1} ,\ldots
$ where $h^0$ corresponds to the trivial link (Fig. 5a) and
$h^i(-A^2-A^{-2})^{i-1}$  corresponds to the link of Fig. 5b. This
notation allows us to embed ${\mathcal S}_{2,\infty}(S^1\times D^2,R)(A)$
in the ring $R[h]$ or equivalently to associate to each link $K \cup L$ in
$S^3$ (as before) a polynomial in $R[h]$. In particular, one has
the following\footnote{Added for e-print:\ Again, if one allows 
the empty link the description would become much simpler.}

\begin{lemma}\label{Lemma 13}
Let $K_1 \cup L_1$ and $K_2 \cup L_2$ be links in $S^3$ with $K_i$
a trivial knot. Then for
 $$K \cup L= (K_1 \cup L_1)
{\#}_{K_i} ( K_2 \cup L_2)\ :
 \ \ L=L_1  \cdot L_2 (-A^2-A^{-2})$$
 (in ${\mathcal S}_{2,\infty}(S^1 \times D^2,R)(A)$).

\end{lemma}
Proof. Notation was chosen so that it holds for initial data. Then
we use the induction.

\begin{lemma} (See \cite{H-P})\label{Lemma 14}
 $A^{-2}$ \parbox{0.9cm}{\psfig{figure=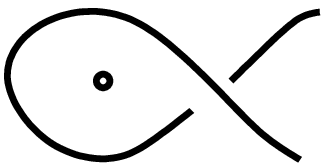,height=0.6cm}}
$+ A^2$ \parbox{0.9cm}{\psfig{figure=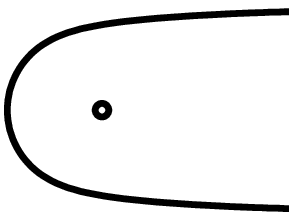,height=0.6cm}}
$=h(A^2 + A^{-2})$
(\parbox{1.3cm}{\psfig{figure=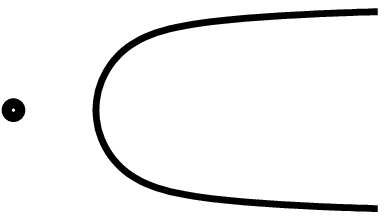,height=0.6cm}})\\
 where a dot $\bullet$ denotes the trivial component
$K$ of a link going vertically to the plane of the projection.
Equivalently we can write $A^{-2}(K \cup L_{++})+ A^2(K \cup
L_{--})=h(A^2+A^{-2})(K \cup L_{00})$ where $(K \cup L_{++}),
(K \cup L_{--})$ and $(K \cup L_{00})$ are link diagrams identical
except the parts shown in Fig. 6. (cf. \cite{H-P}).
\end{lemma}

\ \\
\centerline{\psfig{figure=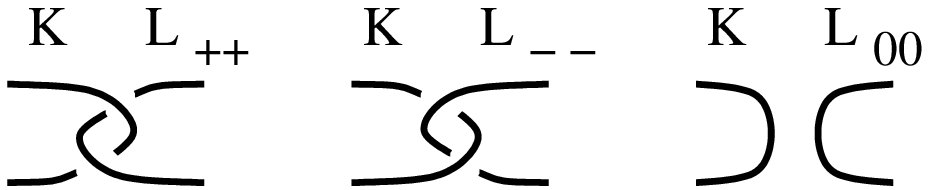,height=2.5cm}}
\begin{center}
Fig. 6
\end{center}

Proof: From the definition
$(-A^{-3})$ \parbox{1.0cm}{\psfig{figure=skeinsmall1.eps,height=0.6cm}} $=$\\
$A$
\parbox{0.8cm}{\psfig{figure=skeinsmall2.eps,height=0.6cm}}
$+A^{-1}$ \parbox{1.2cm}{\psfig{figure=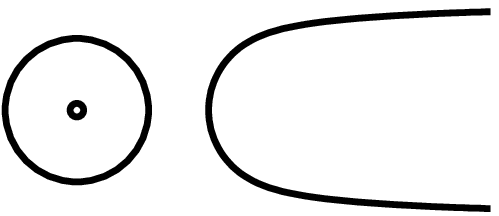,height=0.6cm}} $=A$
\parbox{0.9cm}{\psfig{figure=skeinsmall2.eps,height=0.6cm}}
$+A^{-1}(h(-A^2-A^{-2}))$
(\parbox{1.3cm}{\psfig{figure=skeinsmall3.eps,height=0.6cm}})
\\
and Lemma 14 follows.

It is interesting to analyze the operator of the change of a base
in ${\mathcal S}_{2,\infty}(H_n,R)(A)$. We will consider it in the simplest
case of $H_1=S^1 \times D^2$.\\

Change of a base in ${\mathcal S}_{2,\infty} (S^1 \times D^2, R)\ (A)$.

Let $\mu$ be meridian of $S^1 \times D^2$ and $\lambda$ and
$\lambda'$ its longitudes ($ \lambda'= \lambda+ \mu$) (see Fig. 7).\\

\ \\
\centerline{\psfig{figure=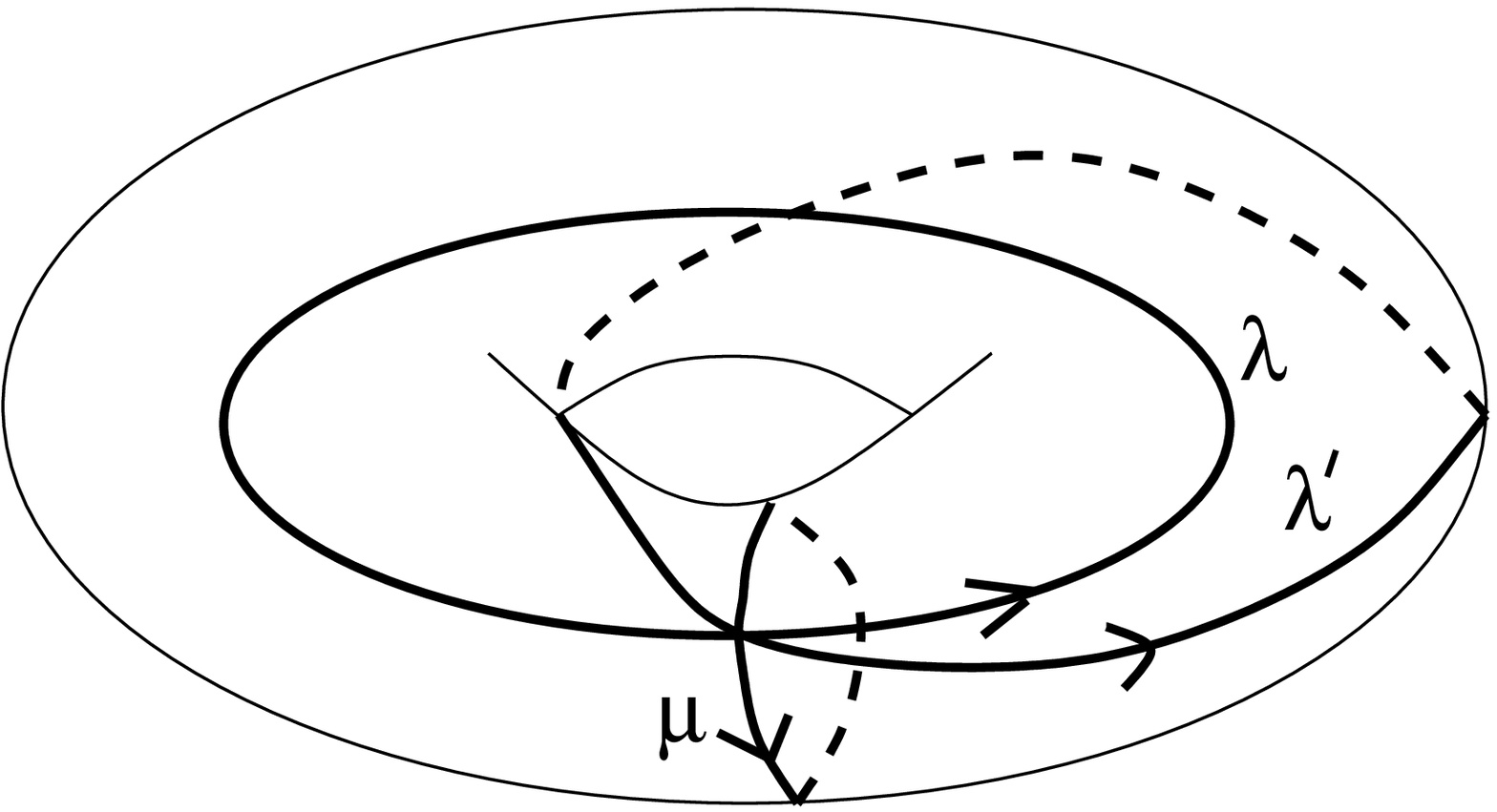,height=5.5cm}}
\begin{center}
Fig. 7
\end{center}
\ \\

Let $1,h, h^2(-A^2-A^{-2}), \ldots $ be a basis of ${\mathcal S}_{2,\infty}
(S^1 \times D^2, R)\ (A)$ with framing $\mu, \lambda $ and $1,
g,g^2(-A^2-A^{-2}), \ldots$ with framing $\mu, \lambda'$.\\

Let $T$ be an automorphism of $R$-module ${\mathcal S}_{2,\infty} (S^1 \times
D^2, R)\ (A)$ such that $T(h^i)=g^i$ in particular $T(1)=1,
T(h)=h$,
$$T(h^2(-A^2-A^{-2}))=A^2h^2(A^2 + A^{-2})+1-A^{-4}$$
$$ T(h^3(-A^2-A^{-2})^2)= A^{2 \cdot 3} h^3(-A^2-A^{-2})+ h(3A^{-2} -2A^6
+2A^{-6})$$

The operator of change of a base in
 ${\mathcal S}_{2,\infty}(S^1 \times D^2, R)\ (A)$
plays an important role when one analyses free $Z_n$ actions on
$S^3$ preserving a link \cite{Pr-2}.

Finally, one should observe that there is the natural surjection from
the skein module onto the (first homology) group algebra:\\

\begin{prop}\label{Proposition 15}
There is the following ring epimorphism from a skein module to a
(first homology) group algebra given by associating to a given
link its homology class\
\begin{itemize}
  \item [(a)]
${\mathcal S}_k(M;R)(r_0, \ldots, r_{k-1})\rightarrow R'[H_1(M;Z)]$
  \item [(b)]${\mathcal S}_k(M, \partial M;R)(r_0, \ldots, r_{k-1})
\rightarrow R'[H_1(M, \partial
  M;Z)]$, where in (a) and (b) $R'$ denotes $R/ (\sum^{k-1}_{i=0}
  r_i)$.
  \item [(c)] ${\mathcal S}_{2,\infty}(M,V_1;R)(A)
\rightarrow R'[H_1(M;Z_2)]$ where $R'$
  denotes $R/(A+A^{-1} -1)$.
\end{itemize}
\end{prop}

\begin{cor}
If $H_1(M,Z)$ is infinite then ${\mathcal S}_k(M)$ is infinitely generated.
\end{cor}

One can build the skein module ${\mathcal S}_{3,\infty}$ using the relation
for the Kauffman \\
polynomial (involving
\parbox{4.8cm}{\psfig{figure=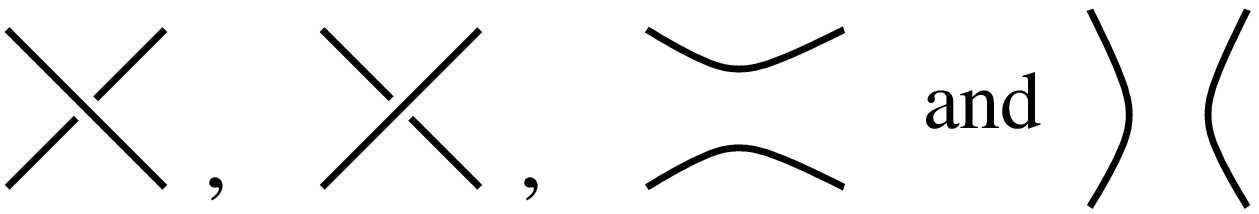,height=0.8cm}})
but we will not analyze it here.


One can consider chromatic skein modules of 3-manifolds,
${\mathcal S}_k^{ch}(M)$. In this case one starts the construction from the
free module over all isotopy classes of chromatic links in $M$ and
each relation involves only one color. The paper \cite{Ho-K} deals with
the dichromatic skein module of $S^3$.

J.~Hoste and the author have also analyzed the homotopy of the
skein module of 3-manifolds (${\mathcal S}^h_k (M,R)(r_0, \ldots,r_{k-1})$)
which is built similarly as the skein module; that is we start
with the set of homotopy classes of oriented links in $M$ (denoted
${\mathcal L}^h(M)$) and we divide the free $R$-module generated
by ${\mathcal L}^h(M)$ (denoted by ${\mathcal L}^h(M,R)$) by the
submodule, $ S^h_{{\mathcal L}(M)}(r_0, \ldots, r_{k-1})$, generated
by linear skein expression $r_0L_0+ \ldots + r_{k-1}L_{k-1}$
assuming that crossings of Fig. 1. involve different components of
$L_{2i}$. That is:
$${\mathcal S}^h_k (M,R)(r_0, \ldots,r_{k-1})=
{\mathcal M}({\mathcal L}^h, R)/S^h_{{\mathcal L}(M)}(r_0, \ldots, r_{k-1})$$
We have been analyzing $S_3^h(M, Z[a^{\mp 1}, z^{\mp 1}])(a^{-1},
-z,a)$ and the computations are easier then in the case of the
skein module. We hope that lens spaces can be distinguished using
this approach.

The first version of this paper was written in May, 1987. In June
1988 we have got the paper \cite{Tu} in which ideas similar to ours are
(independently) developed. Turaev considers also Kauffman modules
($S_{3,\infty}$ in our notation) and computes it for $S^1 \times
D^2$ (the result obtained also in \cite{Ho-K}). We should mention here one
observation from \cite{Tu}. Namely, it is often convenient to consider
skein modules of framed links. In particular ${\mathcal S}_{2, \infty}(M)$
(and ${\mathcal S}_{3,\infty}(M)$) can be defined in that way independently
of the Heegaard decomposition of the manifold $M$.\\

\vspace*{0.1cm}

{\it Added in proof.} Conjecture 7 has been proved by the author in
\cite{Pr-3}.\\

\vspace*{0.3cm}

\noindent {\footnotesize
THE INSTITUTE FOR ADVANCED STUDY, PRINCETON NEW JERSEY 08540 USA}

\end{document}